\documentclass[12pt]{amsart}
\usepackage{bbold}
\usepackage{dsfont}
\usepackage{amsmath}
\usepackage{amsfonts}
\usepackage{amssymb}
\usepackage{graphicx}

\newtheorem{theo}{Theorem}[section]
\theoremstyle{plain}
\newtheorem*{acknowledgement}{Acknowledgment}

\newtheorem{cor}[theo]{Corollary}

\newtheorem{lemma}[theo]{Lemma}

\newtheorem{proposition}[theo]{Proposition}

\begin{document}
\title{A note on the Schur and Phillips lemmas}
\author{Ahmed  Bouziad}
\address[]
{D\'epartement de Math\'ematiques,  Universit\'e de Rouen, UMR CNRS 6085, 
Avenue de l'Universit\'e, BP.12, F76801 Saint-\'Etienne-du-Rouvray, France.}
\email{ahmed.bouziad@univ-rouen.fr}
 
\begin{abstract} It is well-known that every weakly convergent sequence in $\ell_1$
is convergent in the norm topology (Schur's lemma).  Phillips' lemma asserts even more strongly  that
 if a sequence $(\mu_n)_{n\in\mathbb N}$ in $\ell_\infty'$ converges pointwise on $\{0,1\}^\mathbb N$ to $0$,
then its $\ell_1$-projection converges in norm to $0$. In this note we show how   the second category version of
Schur's lemma, for which a short proof is included, can be used  to replace in Phillips' lemma $\{0,1\}^\mathbb N$
 by any of its subsets which contains all finite sets and   having some kind of interpolation property for finite sets.
\end{abstract}

\maketitle
\section{Introduction} 
A classical result of R.S. Phillips \cite{P}, known as  ``Phillips' lemma'', states that for every sequence $(\mu_n)_{n\in\mathbb N}\subset \ell_\infty'$  such that $\lim_n\mu_n(y)=0$ for each $y\in\{0,1\}^\mathbb N$, we have
$\lim_n||\mu_n||_1=0$.  Here $\ell_\infty'$  stands for the conjugate of the Banach space $\ell_\infty$ of all bounded sequences
equipped with the supremum norm  
and $||\mu||_1=\sum_{k\in\mathbb N}|\mu(k)|$ is the $\ell_1$-norm of $\mu\in \ell_\infty'$.  We shall  denote by  $\mu^1$ the $\ell_1$-projection (or the atomic part)
of $\mu$ given by $\mu^1(y)=\sum_{k\in\mathbb N}\mu(k)y(k)$, for $y\in \ell_\infty$. Thus $||\mu||_1$
is the norm  
of $\mu^1$ in  the Banach space $\ell_1$  of absolutely summable sequences; it is well-known
  that $\mu^1\in \ell_1$ and $||\mu^1||\leq ||\mu||$ \cite{Di}.  Schur's lemma \cite{S} 
	can be derived from Phillips' lemma by applying the latter result to sequences $(\mu_n)_{n\in\mathbb N}$ which are already in $\ell_1$, that is, $\mu_n^1=\mu_n$ for every $n\in\mathbb{N}$.
	
	\par
In this note, we shall show by relying on the ``second category version'' of Schur's lemma that it suffices to assume in Phillips' lemma that 
$\lim_n \mu_n (x) = 0 $
for each $ x \in C $, where $ C \subset \{0,1 \}^\mathbb N $ is any set containing finite sets and
having some interpolation property for finite sets (the precise definition is given below).
   Roughly speaking, the key tool for this task is the fact  that
      each $ \mu \in \ell_\infty'$ must   match 
with its $ \ell_1 $-projection on a Baire subspace of the Cantor space $ \{0,1 \}^\mathbb N $. This is expressed in Lemma 3.1 in a form suitable for our purposes. \par
   The second category version  of  Schur's lemma in question asserts that in order that $\lim_n||\mu_n||=0$  (for $(\mu_n)_{n\in\mathbb N}\subset \ell_1$)
it is sufficient that $\lim_n\mu_n(x)=0$ for all $x$ in some non-meager (or second
category) set $C\subset \{0,1\}^\mathbb N$  
  containing all finite sets (see Proposition 2.2). This result is recently proved in \cite[Theorem 3.2]{So}, where
	such a set $C$ is called a Schur family\footnote[2]{Theorem 3.2 in \cite {So} asserts that every non-meager set containing finite sets satisfies Phillips' lemma.  We must admit  that one point pointed out by Bat-Od Battseren in
	Sobota's argumentation eluded us; however,  the proof seems to work  for Schur's lemma.}. We will include here a short 
proof 
of Sobota's theorem; this proof (that we believe is new) is based on a simple ``norming'' property of $\ell_1$ established in Lemma 2.1
below. \par
As usual, $\mathbb N$ stands for the set $\{1,2,\ldots\}$  of positive integers and $\mathbb R$
is the real line; we shall use the standard identification of the power set of $\mathbb N$
with $\{0,1\}^\mathbb N$  equipped  with the product
topology ($\{0,1\}$ being discrete).  We refer to J. Diestal book \cite{Di} for undefined functional-analytic terms.
 \section{Schur's lemma}
Let us say that a set $K\subset \ell_\infty$ is a {\it lower norming set} of $\ell_1$ if for every $x\in \ell_1$
we have $||x||\leq\sup_{u\in K}u(x)$. It is not assumed here  that lower norming sets are subsets of the unit ball of $\ell_\infty$. It is well known and easy to see that $\{-1,1\}^{\mathbb N}$
is a (lower) norming set. We shall say that  $y\in \ell_\infty$ is a {\it null sequence}
if there is $n\in\mathbb N$ such that $y(i)=0$ for every $i\geq n$.\par
We shall give a proof of Schur's lemma based on the following property of $\ell_1$: 
\begin{lemma} Let $T\subset \mathbb R$ be a finite set such
that $T^\mathbb N$ is a lower norming set of $\ell_1$.  Let $C$ be a dense subset of $T^\mathbb N$ {\rm(}with respect
to the product topology{\rm)}.  Then for every nonempty
open subset  $W$ of $C$,  there is a finite set $L\subset \ell_\infty$  of null sequences
such that for each $x\in \ell_1$ we have 
$$||x||\leq \sup_{(u,v)\in W\times L}(u+v)(x),$$
i.e. $W+L$ is lower norming.
\end{lemma}
\begin{proof}  Let us first fix some notations.  For each
$\alpha\in\mathbb R^j$ where $j\in\mathbb N$, let $\pi_\alpha:\ell_\infty\to \ell_\infty$ be the function
defined by $\pi_\alpha(y)(i)=\alpha(i)$ if $i\leq j$ and $\pi_\alpha(y)(i)=y(i)$
for $i>j$. Put $L_\alpha=\{y-\pi_\alpha(y):y\in T^\mathbb N\}$ and note that  $L_\alpha$ is finite and
$z(i)=0$ for each $z\in L_\alpha$ and $i>j$. \par
Now, let $W$ be a nonempty open subset of $T^\mathbb{N}$
 and choose  $\eta>0$, $j\in\mathbb N$, and $\alpha\in \mathbb R^j$
 so that $\cap_{i\leq j}\{y\in T^\mathbb  N:|y(i)-\alpha(i)|<\eta\}\cap C\subset   W$. Let $x\in \ell_1$ and $\varepsilon>0$. There 
 is $y\in T^\mathbb N$  such that $||x||< y(x)+\varepsilon$.
Since the weak* topology on bounded 
subsets of  $\ell_\infty$
 is the topology induced by 
the product topology of  $\mathbb{R}^\mathbb{N}$ and $C$ is dense in $T^\mathbb{N}$,
   there is $y_x\in C$  such that $||x||<y_x(x)+\varepsilon$, that is
 $$||x||< \pi_\alpha(y_x)(x)+(y_x-\pi_\alpha(y_x))(x)+\varepsilon.$$  Since
$C$ is dense  in $T^\mathbb N$,  we have
$\pi_\alpha(y_x)\in \overline W$, which implies that 
 $$||x||\leq \sup_{(u,v)\in W\times L_\alpha}(u+v)(x)+\varepsilon.$$
Since $\varepsilon$ was arbitrary, we get that:
$$\|x\|\le\sup_{(u,v)\in W\times L_\alpha}(u+v)(x).$$
\end{proof}
Let $C_0$ stand for the set of all finite subsets of $\mathbb N$.  The following is established in \cite{So}.
\begin{proposition} {\rm (Schur's Lemma for second category sets)} Let $(\mu_n)_{n\in\mathbb N}\subset \ell_1$ and $C$ a non-meager  set
in $\{0,1\}^\mathbb N$ such that $C_0\subset C$. If $\lim_n\mu_n(x)=0$ for every
$x\in C$, then $\lim_n||\mu_n||=0.$
\end{proposition}
 \begin{proof}   Let $\varepsilon>0$. To simplify, replacing
$C$ by $\{x-y: x\in C, y\in C_0\}$, we assume that $C$
is a dense second category subspace of $\{-1,0,1\}^\mathbb N$. Since $C\subset \cup_{n\in\mathbb N}\cap_{k\geq n}\{y\in C: |y(\mu_k)|\leq \varepsilon\}$,
there is a nonempty open subset $W$ of $C$
and 
  $n_0>0$ such that  $|y(\mu_n)|\leq\varepsilon$ for every
 $n\geq n_0$ and $y\in W$.  
By Lemma 2.1, there is a finite  $L\subset \ell_\infty$ of null sequences  so that 
$$||\mu_n||\leq \sup_{y\in W}|y(\mu_n)|+\sum_{y\in L}|y(\mu_n)|$$
for each $n\in\mathbb N$. Since $L$ is a finite set of null sequences and $C$
contains all finite sets, we have $\lim_n y(\mu_n)=0$ for every $y\in L$,  thus $\limsup_n||\mu_n||\leq \varepsilon$.
\end{proof}
\section{Phillips' lemma}
 For  a given sequence $\sigma=(\mu_n)_{n\in\mathbb N}$  in $\ell_\infty'$, 
 let $\mathcal G_\sigma$ stand  for the set of all $A\subset\mathbb N$ such that $\mu_n(E)=\mu_n^1(E)$ for every
$E\subset A$ and $n\in\mathbb N$. It is easy to see that $\mathcal G_\sigma$ contains $C_0$ and is closed
under finite union and subsets, that is,  $\mathcal G_\sigma$ is an ideal on $\mathbb N$.  \par
In what follows, we shall write $A\subset^* B$ to say that $A$
is almost contained in $B$, that is, $A\setminus B$ is finite. Recall
that a space $X$ is said to be a {\it Baire space} if every intersection of countably many dense open sets in $X$ is dense in $X$. The key observation mentioned in the introduction  is the following  fact: 
 
 \begin{lemma} Let 
$C\subset\{0,1\}^\mathbb N$. Suppose that for every sequence
$(F_n)_{n\in\mathbb N}$ of finite disjoint subsets of $\mathbb N$, there exist $A\in C$ and  an infinite set $I\subset \mathbb N$
such that $\cup_{i\in I}F_n\subset A\subset^*\cup_{n\in \mathbb N}F_n$. Then, for every $\sigma=(\mu_n)_{n\in\mathbb{N}}\subset\ell_\infty'$,
  $C\cap \mathcal G_\sigma$  is a non-meager subset of
$\{0,1\}^\mathbb N$. In particular, $\mathcal G_\sigma$
is a Baire space.
\end{lemma}
 \begin{proof} Let $(\mu_n)_{n\in\mathbb{N}}\subset\ell_\infty'$. Since every $\mu_n^1$ is $\sigma$-additive, $\mathcal G_\sigma$ is none other than the set of $A\subset\mathbb N$
such that each $\mu_n$
is $\sigma$-additive on the $\sigma$-algebra of all subsets of $A$. Let $(E_n)_{n\in\mathbb N}$
be finite intervals partition of $\mathbb N$, $X\subset\mathbb N$ and let us show that there is $A\in C\cap\mathcal G_\sigma$ 
such that  $E_n\cap X=E_n\cap A$ for infinitely many  $n\in \mathbb N$; this will imply that $C\cap \mathcal G_\sigma$
is nonemeager (see Theorem 5.2 in \cite{B}). Put $F_n=X\cap E_n$, $n\in\mathbb N$.\par
Since each $\mu_n$
is exhaustive, i.e. $\lim_k\mu_n(D_k) = 0$
for every disjoint sequence $(D_k)_{k\in\mathbb N}$ in $\{0,1\}^\mathbb N$ (\cite{Di}), there is a decreasing sequence $(I_n)_{n\in\mathbb N}$ of infinite sets such that 
for every $l,k\in\mathbb{N}$, $l\le k$, we have $|\mu_l(E)|<1/k$
for  each $E\subset\cup_{n\in I_k}E_n$. Let $J\subset \mathbb N$ be an infinite set such
that $J\subset^*I_n$ for every $n\in\mathbb N$. There exist  $A\in C$ 
and an infinite set $I\subset J$  such that $\cup_{n\in I}F_n\subset A\subset^*\cup_{n\in J}F_n$. Then, taking $k\in \mathbb N$
so that $A\setminus\cup_{n\in J}F_n\subset\cup_{i\leq k}E_k$, we get
$E_n\cap A=F_n$ for every $n\in I$ such that $n> k$.  \par
To conclude, let us check that  $\cup_{n\in J}E_n$ is  in 
$\mathcal G_\sigma$, this will imply that $A$ is also in $\mathcal G_\sigma$. We will show that each $\mu_n$ is $\sigma$-additive on the $\sigma$-algebra of subsets of $\bigcup_{m\in J}E_m$. Let $l\in\mathbb N$.
Let  $(A_n)_{n\in\mathbb N}\subset \cup_{n\in J}E_n$  be an arbitrary  decreasing sequence such that $\cap_{n\in\mathbb N}A_n=\emptyset$
and let us show that $\lim_n\mu_l(A_n)=0$; since $\mu_l$
is bounded and finitely additive, it will follow  that $\mu_l$
is $\sigma$-additive on the power set of $\bigcup_{m\in J}E_m$.  Let $\varepsilon>0$.  Choose
$k_0\in\mathbb N$ so that $l\leq k_0$ and  $1/k_0<\varepsilon$, and let
 $m_0\in\mathbb N$ be such that $J\subset\{0,\ldots,m_0\} \cup I_{k_0}$. Since $\cap_{n\in\mathbb N}A_n=\emptyset$ and $\cup_{n\leq m_0}E_n$
is finite,
there is $n_0\in\mathbb N$
such that $A_{n_0}\cap(\cup_{n\leq m_0}E_n)=\emptyset$. Then $A_{n_0}\subset \cup_{n\in I_{k_0}}E_n$; for
if $x\in A_{n_0}$, then $x\in E_i$ for some $i\in J$ and, since $x\not\in \cup_{n\leq m_0}E_n$, we have $i\in I_{k_0}$, and thus
 $x\in \cup_{n\in I_{k_0}}E_{n}$. It follows that
$|\mu_l(A_n)|<1/k_0$ for every $n\geq n_0$. 
\end{proof}

Let us emphasize  that the essential point in the above proof shows  that  for any sequence
$(F)_{n\in\mathbb N}$ of disjoint (finite) subsets of $\mathbb N$, and every $\mu\in \ell_\infty'$, there is an infinite set $I\subset \mathbb N$ such 
that $\mu$ is $\sigma$-additive on  the $\sigma$-algebra of all subsets of $\cup_{n\in I}F_n$; this is much stronger than
saying that
$\mu$ is $\sigma$-additive on  the $\sigma$-algebra generated by $\{F_n:n\in I\}$ (\cite{D}).\par
Let us say that a set  $C\subset\{0,1\}^\mathbb N$ is an  {\it ipfset} (for  ``interpolation property for finite sets'')  if
for any sequence $(F_n)_{n\in\mathbb N}\subset C_0$ of disjoint sets, there are an infinite set $I\subset\mathbb N$
and $A\in C$ such that $\cup_{n\in I}F_n\subset A\subset^*\cup_{n\in\mathbb N}F_n$. It goes without saying that the introduction of the ipfset property
is motivated here by the validity of
  Lemma 3.1.  The ipfset property resembles without being strictly comparable to some other properties (of  interpolation or  completeness type concerning rings) that can be found in the literature; for instance, see 
\cite{C,H1,H,W} and the references therein (this list is far from being exhaustive). 
 \begin{proposition} {\rm (Phillips' lemma for ipfsets)} Let $\sigma=(\mu_n)_{n\in\mathbb N}$ be a sequence in $\ell_\infty'$
 and $C\subset \{0,1\}^\mathbb N$ be an ipfset containing $C_0$  such that
  $\lim\mu_n(x)=0$ for each $x\in C$. Then
 $\lim_n||\mu_n^1||=0$.
  \end{proposition}
 \begin{proof} Since by  Lemma  3.1, the subset $C\cap \mathcal G_\sigma$ of $\{0,1\}^\mathbb N$
is of the second category,   Schur's Lemma (Proposition 2.2) allows  to conclude that $\lim_n||\mu_n^1||=0$.
 \end{proof} 
One can easily see that 
every second category subspace of $\{0,1\}^\mathbb N$ containing $ C_0$ which is either  closed by taking subsets
or totally non-meager (i.e. all its  closed subspaces are Baire spaces)  is an ipfset.  Therefore, we get
 from Proposition 3.2:
\begin{cor}  Let $\sigma=(\mu_n)_{n\in\mathbb N}$ be a sequence in $\ell_\infty'$
 and $C\subset \{0,1\}^\mathbb N$  containing $C_0$  such that
  $\lim\mu_n(x)=0$ for each $x\in C$. If $ C $ is either totally non-meager
or second category and closed by taking   subsets, then
 $\lim_n||\mu_n^1||=0$.
\end{cor}

A classical theorem due to Sierpi\`nski asserts that every ultrafilter  $\mathcal F $ on $\mathbb N$
is a Baire subspace of $\{0,1\}^\mathbb N$. Let us note, to  conclude, that this result follows from Lemma 3.1. Indeed,
  let  $\sigma =(\mu_n)_{n\in\mathbb N}\subset \ell_\infty'$ be the constant sequence where, for every $n\in\mathbb N$
	and $x\in \ell_\infty$,
$\mu_n(x)$ is the limit  along  $\mathcal F$ of the sequence $(x(k))_{k\in\mathbb N}$. Then 
$\mathcal G_\sigma$ is the co-ideal $\mathcal I$ of $\mathcal F$; therefore, by Lemma 3.1, $\mathcal I$ (equivalently, $\mathcal F$) is a Baire space.
\begin{acknowledgement} {\rm The author thanks the referee for valuable comments and suggestions which improved the presentation of this manuscript.}
 \end{acknowledgement}


\begin{thebibliography}{99}
\bibitem{B} A. Blass,{\it  Combinatorial cardinal characteristics of the continuum}, Handbook of set theory. Vols. 1, 2, 3, Springer, Dordrecht, 2010, pp. 395--489.
\bibitem{C} C. Constantinescu, {\it
On Nikod\'ym's boundedness theorem,}
Libertas Math. 1 (1981), 51–-73. 
\bibitem{Di} J. Diestel,
Sequences and series in Banach spaces,
Graduate Texts in Mathematics, 92. Springer-Verlag, New York, 1984.
\bibitem{D} L. Drewnowski,{\it 
Uniform boundedness principle for finitely additive vector measures} (Russian summary),
Bull. Acad. Polon. Sci. S\'er. Sci. Math. Astronom. Phys. 21 (1973), 115--118. 
\bibitem{H1} R.G. Haydon, {\it  
Boolean rings that are Baire spaces},
Serdica Math. J. 27 (2001), no. 2, 91--106. 

\bibitem{H} R. G.  Haydon, 
{\it A nonreflexive Grothendieck space that does not contain $\ell_\infty$},
Israel J. Math. 40 (1981), no. 1, 65--73. 
\bibitem{P}R. S. Phillips, {\it On linear transformations}, Trans. Amer. Math. Soc. 48 (1940), 516--541.

\bibitem{So} D. Sobota, {\it 
On families of subsets of natural numbers deciding the norm convergence in $\ell_1$},
Proc. Amer. Math. Soc. 146 (2018), no. 4, 1673--1680. 

\bibitem{S} J. Schur, {\it \"Uber lineare Transformationen in der Theorie der unendlichen Reihen }(German), 
J. Reine Angew. Math. 151 (1921), 79--111.
\bibitem{W} B.B. Wells, Jr.,{\it 
Weak compactness of measures,}
Proc. Amer. Math. Soc. 20 (1969), 124--130.
\end{thebibliography}
\end{document}